\documentclass[preprint,12pt]{article}
\usepackage{amssymb}
\usepackage{amsmath}
\usepackage{pdfpages}
\usepackage{amsfonts}
\usepackage{float}
\usepackage{fullpage}
\usepackage{xcolor}
\usepackage{cleveref}

\newcommand{\Dcal}{\mathcal{D}}

\title{Filtered Neural Galerkin model reduction schemes for efficient propagation of initial condition uncertainties in digital twins}
\author{Zhiyang Ning\thanks{Courant Institute of Mathematical Sciences, New York University} \and Benjamin Peherstorfer\footnotemark[1]}

\begin{document}

\maketitle

\begin{abstract}
Uncertainty quantification in digital twins is critical to enable reliable and credible predictions beyond available data. 
A key challenge is that ensemble-based approaches can become prohibitively expensive when embedded in control and data assimilation loops in digital twins, even when reduced models are used. 
We introduce a reduced modeling approach that advances in time the mean and covariance of the reduced solution distribution induced by the initial condition uncertainties, which eliminates the need to maintain and propagate a costly ensemble of reduced solutions. 
The mean and covariance dynamics are obtained as a moment closure from Neural Galerkin schemes on pre-trained neural networks, which can be interpreted as filtered Neural Galerkin dynamics analogous to Gaussian filtering and the extended Kalman filter.
Numerical experiments demonstrate that filtered Neural Galerkin schemes achieve more than one order of magnitude speedup compared to ensemble-based uncertainty propagation. 
\end{abstract}

\section{Introduction}\label{sec:Intro}
Uncertainty quantification in digital twins is critical to predict reliably and credibly beyond available data  \cite{national2024foundational,kapteyn2021probabilistic,TORZONI2024116584,pash2025predictivedigitaltwinsquantified}.  
In this work, we focus on propagating initial condition uncertainties through a model that describes part of the physical asset of a digital twin. 
Specifically, we consider initial condition uncertainty propagation through nonlinear reduced models, which are important building blocks to enable control and data assimilation in digital twins \cite{Rozza2008,doi:10.1137/130932715,interpbook,P22AMS,annurev:/content/journals/10.1146/annurev-fluid-121021-025220}. 
Taking an ensemble of reduced-model solutions to represent the uncertainties seems at first tractable because reduced models incur low computational costs per solve; however, when wrapped into control and data assimilation loops in digital twins, even ensembles of reduced solutions can become computationally prohibitive. For example, consider taking ten ensemble members, then the costs increase already by one order of magnitude. In our experiments, even larger ensemble sizes are needed to accurately estimate mean and covariance from an ensemble to quantify uncertainties. 
In contrast, in this work, we propose a filtered reduced modeling approach that propagates forward the mean and covariance of the distribution of the reduced solutions induced by the initial condition uncertainties. 
We develop filtered reduced models based on Neural Galerkin schemes \cite{BRUNA2024112588,BP23RSNG,WEN2024134129,JSP24Handbook} with pre-trained neural networks \cite{pmlr-v235-berman24b} that lead to nonlinear approximations and so are effective and efficient in quantifying uncertainties in  transport- and wave-dominated problems \cite{P22AMS}. 
We call these the filtered Neural Galerkin equations by analogy with Gaussian filtering and the extended Kalman filter: the filtered equations are the first-order---Gaussian---moment closure of the Neural Galerkin dynamics. 
The computational costs per time step increase by a factor that scales only with the reduced dimension, which leads to more than one order of magnitude cost reduction compared to ensemble-based uncertainty quantification with reduced models in our experiments. 

There is a range of methods for propagating uncertainties in initial conditions through computational models \cite{StochCompBook,HandbookUQ} as well as using reduced models for uncertainty quantification  \cite{doi:10.1137/120876745,doi:10.1137/151004550,doi:10.1137/16M1082469}. 
There are ensemble-based methods such as Monte Carlo and collocation methods, which can be combined with multi-level and multi-fidelity concepts to reduce computational costs \cite{Giles_2015,doi:10.1137/16M1082469}. Instead, we specifically want to propagate uncertainties through reduced models to account for errors in reduced initial-condition approximations. 
Another line of work learns transport maps to push forward a source distribution to a target distribution \cite{Marzouk2016}. In our case of propagating initial condition uncertainties, one could learn a transport map from the initial condition to the distribution induced by the flow through the reduced model. However, this requires a separate training procedure as well as training data to construct the map. 
Closer to our work is \cite{SAPSIS20092347,doi:10.1137/16M1109394,UECKERMANN2013272} that introduces reduced models based on dynamic orthogonal decomposition with dynamically changing reduced basis for stochastic problems by deriving moment equations. Instead, we consider pre-trained nonlinear parametrizations that are not updated online. In \cite{HILLIARD2025113994,NGDA2025}, the authors combine the same dynamics as used for Neural Galerkin schemes with data assimilation; however, there the goal is achieving a higher approximation accuracy instead of uncertainty propagation as in our case.  
There is a large body of literature on Kalman filtering and reduced modeling. First, there are reduced Kalman filters  \cite{NAGPAL01061987,StateEstimationUsingaReducedOrderKalmanFilter,ROMUCF,4434882,HaasdonkRBKF} with the goal of reducing the costs of Kalman filtering for problems with high dimensional state dimensions.  While we derive the filtered Neural Galerkin equations with analogous steps as used to derive the extended Kalman filter, we skip the update step that assimilates observations and we operate already on the reduced state. 
Another large body of work develops low-rank solvers for approximately solving matrix equations that arise in control and filtering  \cite{https://doi.org/10.1002/nla.622,doi:10.1137/09075041X,1272745}. In contrast, the dynamics we want to filter are already low dimensional but not closed and we filter them to close them rather than to reduce costs. 
There is work that combines model reduction with variational formulations of data assimilation \cite{https://doi.org/10.1002/nme.4747,Maday2013} and the ensemble Kalman filter \cite{doi:10.1137/16M1078598,Silva2023,doi:10.1137/24M1653690}; however, we are not considering data assimilation because we skip the update step that assimilates observations with states.  
We note that there is also a line of work in reduced modeling that aims to match moments of transfer functions of dynamical systems \cite{MR1850408,doi:10.1137/130932715}, but these methods target frequency-domain approximations.
The work \cite{GUO2022115336} introduces a Bayesian approach to learning reduced models from data, which can propagate uncertainties over time. In contrast, we have a reduced model given via Neural Galerkin schemes on pre-trained neural networks and want to propagate initial condition uncertainties without learning a new model.

We derive Neural Galerkin schemes that propagate forward uncertainties in the initial condition. One source of uncertainty can be partial knowledge of the initial condition, in which case the initial condition is modeled as a distribution to account for the partial knowledge. Another source of uncertainty is specific to reduced modeling and stems from approximating the initial condition in a reduced representation that introduces errors. The starting point for us is that a distribution is given that represents such initial condition uncertainties. 
To derive the filtered dynamics to propagate forward initial condition uncertainties, we follow similar steps as in the extended Kalman filter: We derive the probability flow of the parameters induced by the Neural Galerkin dynamics and the distribution that represents the initial condition. We then derive the equations for the mean and covariance---first and second moments---and close them with a Gaussian filter that corresponds to linearizing the Neural Galerkin dynamics over time. The resulting filtered Neural Galerkin equations can be integrated in time with standard time integration schemes. 
We show that the costs per time step increase by a factor that scales with the reduced dimension. The low cost complexity is demonstrated with numerical experiments where the filtered Neural Galerkin schemes achieve more than one order of magnitude speedup compared to ensemble-based reduced methods.

This manuscript is structured as follows. In \Cref{sec:Prelim}, we provide preliminaries about Neural Galerkin schemes and pre-training neural networks for model reduction and give a problem formulation. The filtered Neural Galerkin schemes are introduced in \Cref{sec:FNG} and numerical experiments are presented in \Cref{sec:NumExp}. \Cref{sec:Conc} provides concluding remarks.

\section{Preliminaries}\label{sec:Prelim}
We discuss preliminaries that cover Neural Galerkin schemes \cite{BRUNA2024112588,BP23RSNG,WEN2024134129} and pre-trained neural networks for model reduction \cite{pmlr-v235-berman24b}. We also provide a problem formulation. 

\subsection{Parametrized evolution equations}
Let us consider a time-dependent partial differential equation (PDE)
\begin{equation}\label{eq:Prelim:PDEProblem}
\begin{aligned}
\partial_t q(t, x; \mu) = f(x, q; \mu)\,,\\
q(0, x; \mu) = q_0(x; \mu)\,,
\end{aligned}
\end{equation}
that depends on a parameter $\mu \in \mathcal{Q} \subseteq \mathbb{R}^{d'}$. We denote the solution field as $q: [0, T] \times \Omega \times \mathcal{Q} \to \mathbb{R}$, which evolves over time $t \in [0, T]$ and depends on the spatial coordinate $x \in \Omega \subseteq \mathbb{R}^d$ and the parameter $\mu$. 
The initial condition is $q_0: \Omega \times \mathcal{Q} \to \mathbb{R}$. 
The right-hand side function $f$ can contain partial derivatives of $q$ with respect to the spatial coordinate.

\subsection{Neural Galerkin schemes}\label{sec:Prelim:NG}
Let $\hat{q}: \mathbb{R}^p \times \Omega \to \mathbb{R}$ be a function that depends on a $p$-dimensional weight vector $\theta(t, \mu) \in \mathbb{R}^p$, which we want to use as a finite-dimensional parametrization of solutions $q$ of the PDE problem \eqref{eq:Prelim:PDEProblem}. 
Notice that the weight vector $\theta(t, \mu)$ depends on time $t$ and parameter $\mu$. 
In the following, we are interested in functions $\hat{q}$ that have a nonlinear dependence on the weight vector $\theta(t, \mu)$ such as neural networks. 
Invoking Neural Galerkin schemes introduced in \cite{BRUNA2024112588} and further developed in \cite{BP23RSNG,WEN2024134129}, we formulate the residual function $r: \mathbb{R}^p \times \mathbb{R}^p \times \Omega \times \mathcal{Q} \to \mathbb{R}$ as
\begin{equation}\label{eq:Prelim:ResidualFun}
r(\theta(t, \mu), \dot{\theta}(t, \mu), x; \mu) = \nabla_{\theta}\hat{q}(\theta(t, \mu), x)^{\top}\dot{\theta}(t, \mu) - f(x, \hat{q}(\theta(t, \mu), \cdot); \mu)\,,
\end{equation}
where $\dot{\theta}(t, \mu) \in \mathbb{R}^{p}$ is the time derivative of $\theta(t, \mu)$. 
In Neural Galerkin schemes, the time derivative $\dot{\theta}(t, \mu)$ is determined via the Dirac-Frenkel variational principle \cite{dirac_1930,Frenkel1934,Lubich2008}\,,
\begin{equation}\label{eq:Prelim:NGOrthogonalityCondition}
\langle \partial_{\theta_i}\hat{q}(\theta(t, \mu), \cdot), r(\theta(t, \mu), \dot{\theta}(t, \mu), \cdot; \mu)\rangle = 0\,,\qquad i = 1, \dots, p\,,
\end{equation}
which sets the residual \eqref{eq:Prelim:ResidualFun} orthogonal in the $L^2$ inner product $\langle \cdot, \cdot \rangle$ to the test space spanned by the components of the gradient 
\[
\nabla_{\theta}\hat{q}(\theta, \cdot) = [\partial_{\theta_1}\hat{q}(\theta, \cdot), \dots, \partial_{\theta_p}\hat{q}(\theta, \cdot)]^{\top}\,.
\]
We refer to \cite{JSP24Handbook} for a survey on Neural Galerkin schemes and to \cite{doi:10.1137/21M1415972,PhysRevE.104.045303,24OTDDTO} for other techniques related to Neural Galerkin schemes. 
System \eqref{eq:Prelim:NGOrthogonalityCondition} can be rewritten as a dynamical system in $\theta(t, \mu)$,   
\begin{equation}\label{eq:Prelim:NGODE}
\frac{\mathrm d}{\mathrm dt} \theta(t, \mu) = v(\theta(t, \mu); \mu)\,, 
\end{equation}
with the velocity field
\begin{equation}\label{eq:Prelim:VField}
v(\theta(t, \mu); \mu) = M(\theta(t, \mu))^{-1}F(\theta(t, \mu); \mu)\,,
\end{equation}
where the components of the matrix $M(\theta(t, \mu)) \in \mathbb{R}^{p \times p}$ and the vector $F(\theta(t, \mu); \mu) \in \mathbb{R}^p$ are 
\[
\begin{aligned}
M_{ij}(\theta) = & \langle \partial_{\theta_i} \hat{q}(\theta, \cdot), \partial_{\theta_j} \hat{q}(\theta, \cdot) \rangle\,,\qquad i, j = 1, \dots, p\,,\\
F_i(\theta; \mu) = &\langle \partial_{\theta_i}\hat{q}(\theta, \cdot), f(\cdot, \hat{q}(\theta, \cdot); \mu)\rangle\,,\qquad i = 1, \dots, p\,.
\end{aligned}
\]
Notice that an evaluation of the velocity field $v$ defined in \eqref{eq:Prelim:VField} incurs a linear solve with $M(\theta(t, \mu))$. In fact, the evaluation of the velocity field is the solution to a least-squares problem for which \eqref{eq:Prelim:VField} corresponds to the normal equations. 
For a sufficiently regular velocity field $v$ (e.g., Lipschitz in $\theta$), system \eqref{eq:Prelim:NGODE} is an ordinary differential equation with the flow map $\Phi_t^{(\mu)}: \mathbb{R}^p \to \mathbb{R}^p$. The flow map $\Phi_t^{(\mu)}$ maps $\theta(0, \mu)$ at time $t = 0$ to $\Phi_t^{(\mu)}(\theta(0, \mu)) = \theta(t, \mu)$ at time $t$.
To obtain the initial weight vector $\theta(0, \mu)$, the initial condition $q_0$ is fitted over a set $\{x_1, \dots, x_M\} \subset \Omega$ of collocation points with the mean-squared error.

\subsection{Pre-training neural networks with CoLoRA layers}
In the following, we will use Neural Galerkin schemes with continuous low-rank adaptation (CoLoRA) neural networks introduced in \cite{pmlr-v235-berman24b} that can be pre-trained on data. 
Motivated by concepts in model reduction and surrogate modeling \cite{Rozza2008,doi:10.1137/130932715,interpbook,P22AMS,annurev:/content/journals/10.1146/annurev-fluid-121021-025220},  trajectory data are collected over training parameters $\mu_1, \dots, \mu_m \in \mathcal{Q}$ from a high-fidelity numerical model of \eqref{eq:Prelim:PDEProblem} with a one-time high cost. The training data are then used to pre-train the CoLoRA neural networks in an offline phase so that only a small number of weights need to be computed with Neural Galerkin schemes in an online phase when one seeks an approximation at a new parameter $\mu \in \mathcal{Q}$. 
Following \cite{pmlr-v235-berman24b}, the training data are given in form of $m$ trajectories
\begin{equation}\label{eq:Prelim:TrainingTraj}
q(\cdot, \cdot; \mu_1), \dots, q(\cdot, \cdot; \mu_m): [0, T] \times \Omega \to \mathbb{R}\,,
\end{equation}
corresponding to the $m$ training parameters $\mu_1, \dots, \mu_m \in \mathcal{Q}$. The training data are typically obtained from numerical simulations of the PDE problem \eqref{eq:Prelim:PDEProblem}. Because we expect them to be of high fidelity, our notation does not distinguish between the PDE solution in, e.g., a variational sense, and the numerical solution. 
The functions given in  \eqref{eq:Prelim:TrainingTraj} represent training trajectories that can be evaluated at times $t \in [0, T]$ and coordinates $x \in \Omega$. For a set 
\begin{equation}
\Omega_{\text{train}} \subset \Omega\label{eq:Prelim:OmegaTrain}
\end{equation}
with a finite number of elements, it will be convenient to define the sets
\begin{equation}\label{eq:Dtmu}
\mathcal{D}(t, \mu) = \{(x, q(t, x; \mu)) \,|\, x \in \Omega_{\text{train}}\}\,.
\end{equation}
If $t = 0$, then the set $\mathcal{D}(0, \mu) = \{(x, q_0(x; \mu)) \,|\, x \in \Omega_{\text{train}}\}$ contains evaluations of the initial condition $q_0$. 
A CoLoRA network depends on the weight vector
\begin{equation}\label{eq:Prelim:CoLoRADecomposition}
\theta_{\text{total}}(t, \mu) = [\theta_{\text{off}}; \theta(t, \mu)]\,,
\end{equation}
which is decomposed into an offline weight vector $\theta_{\text{off}} \in \mathbb{R}^n$ that is independent of time $t$ and parameter $\mu$ and an online weight vector $\theta(t, \mu) \in \mathbb{R}^p$ that changes with time and parameter. The offline weight vector $\theta_{\text{off}}$ is learned during pre-training and then fixed online when the online weight vector $\theta(t, \mu)$ is computed with Neural Galerkin schemes. 
The decomposition \eqref{eq:Prelim:CoLoRADecomposition} is induced by the CoLoRA network architecture, which builds on CoLoRA layers introduced in \cite{pmlr-v235-berman24b} as
\begin{equation}\label{eq:Prelim:CoLoRALayer}
C_{\ell}(y) = W_{\ell}y + \theta_{\ell}(t, \mu)A_{\ell}B_{\ell}y + b_{\ell}\,,
\end{equation}
where $W_{\ell}$ is a matrix and $b_{\ell}$ is a vector of appropriate size. The index $\ell \in \mathbb{N}$ stands for the index of the layer. The matrix given by the product $A_{\ell}B_{\ell}$ is of rank $r$ with $A_{\ell}$ having $r$ columns and $B_{\ell}$ having $r$ rows. The coefficient $\theta_{\ell}(t, \mu) \in \mathbb{R}$ is a scalar. Typically the rank $r$ is small compared to the size of the matrix $W_{\ell}$ and the vector $b_{\ell}$. We stress that other variants of CoLoRA layers are introduced in \cite{pmlr-v235-berman24b} but we will focus on CoLoRA layers of the form given in \eqref{eq:Prelim:CoLoRALayer}. 

A CoLoRA network $ \hat{q}(\theta(t, \mu), \cdot; \theta_{\text{off}}): \Omega \to \mathbb{R}$ composes $\ell = 1, \dots, L$ CoLoRA layers,
\[
\hat{q}(\theta(t, \mu), x; \theta_{\text{off}}) = C_{\text{out}}(\sigma(C_L( \sigma( C_{L-1} ( \dots \sigma (C_1(x)) \dots )))))\,,
\]
where $\sigma$ is an activation function. 
The function $C_{\text{out}}$ is the output layer of the form 
\begin{equation}\label{eq:Prelim:CoLoRAOutputLayer}
C_{\text{out}}(y) = w_{L+1}^{\top}y + \theta_{L+1}(t, \mu) A_{L+1}B_{L+1}y + b_{L + 1}\,,
\end{equation}
which also depends on an online weight $\theta_{L + 1}(t, \mu)$ where $A_{L+1}$ is a scalar and $B_{L+1}$ is a row vector. The weight $w_{L+1}$ is a vector. 
The online weight vector $\theta(t, \mu) = [\theta_1(t, \mu), \dots, \theta_L(t, \mu), \theta_{L + 1}(t, \mu)]^{\top}$ has dimension $p = L + 1$. 

A CoLoRA network is pre-trained on the training data \eqref{eq:Prelim:TrainingTraj}. The pre-training is achieved via a hyper-network: Let $h_{\psi}: [0, T] \times \mathcal{Q} \to \mathbb{R}^p$ be a neural network with weights $\psi \in \mathbb{R}^{p'}$. Typically the hyper-network is a fully connected feedforward network with a few layers only. The hyper-network is only used for the pre-training; see \cite{pmlr-v235-berman24b} for details. 
Let us now consider the objective
\begin{equation}\label{eq:Prelim:CoLoRAObj}
\mathcal{L}(\theta_{\text{off}}, \psi) = \sum_{i = 1}^m \sum_{\substack{x \in \Omega_{\text{train}}\\t \in \{t_0, \dots, t_K\}}} \frac{|q(t, x; \mu_i) - \hat{q}(h_{\psi}(t, \mu_i), x; \theta_{\text{off}})|^2}{|q(t, x; \mu_i)|^2}\,,
\end{equation}
where the set $\Omega_{\text{train}} \subset \Omega$ is the finite subset of $\Omega$ defined in \eqref{eq:Prelim:OmegaTrain} and $0 = t_0 < t_1 < \dots < t_K = T$ are discrete time steps in the time interval $[0, T]$. Notice that the evaluations of $q$ used in the objective $\mathcal{L}$ are given by the training data \eqref{eq:Prelim:TrainingTraj} corresponding to $\mu_1, \dots, \mu_m$. 
Once the CoLoRA network $\hat{q}$ is pre-trained with the objective \eqref{eq:Prelim:CoLoRAObj} on the data \eqref{eq:Prelim:TrainingTraj}, only the online weight vector $\theta(t, \mu)$ needs to be computed to approximate PDE solutions at a new parameter $\mu$ over time $[0, T]$. 
For example, in \cite{pmlr-v235-berman24b}, Neural Galerkin schemes, as described in \Cref{sec:Prelim:NG}, are used to determine the online weight vector $\theta(t, \mu)$.

\subsection{Problem formulation}
We are interested in the scenario that the initial condition is a random function that is represented by $\hat{q}(\Theta(0, \mu), \cdot): \Omega \to \mathbb{R}$ with a random variable $\Theta(0, \mu)$. The distribution of the random variable $\Theta(0, \mu)$ is $\pi_0^{(\mu)}$, which is supported on $\mathbb{R}^p$. 
Such a situation arises, for example, if the initial condition $q_0$ is modeled stochastically due to partial knowledge and the distribution $\pi_0^{(\mu)}$ is fitted so that $\hat{q}(\Theta(0, \mu), \cdot)$ is close in a statistical sense to $q_0$ for $\Theta(0, \mu) \sim \pi_0^{(\mu)}$. 
We also find this situation when $q_0$ is deterministic but a distribution $\pi_0^{(\mu)}$ of weights is fitted so that the distribution accounts for uncertainties in the fitting process, e.g., when fitting Bayesian neural networks to data \cite{Blei03042017}. Our goal is predicting how the distribution $\pi_0^{(\mu)}$ is propagated as $\pi_t^{(\mu)}$ over time $t$ when following the Neural Galerkin dynamics \eqref{eq:Prelim:NGODE}. 
The random variable $\Theta(t, \mu) \sim \pi_t^{(\mu)}$ and the function $\hat{q}(\Theta(t, \mu), \cdot)$ should represent how the initial condition uncertainties propagate over time.  
We specifically want to avoid ensemble-based approaches because computing ensembles of  solutions can be computationally expensive; see the discussion in \Cref{sec:Intro}.

\section{Filtered Neural Galerkin schemes}\label{sec:FNG}
We introduce filtered Neural Galerkin schemes that evolve the mean and covariance of the weight distribution $\pi_t^{(\mu)}$ induced by Neural Galerkin dynamics with an initial distribution $\pi_0^{(\mu)}$.
The filtered equations are derived analogously to the extended Kalman filter \cite{Sarkka_2013,Sanz-Alonso_Stuart_Taeb_2023} by first considering the continuity equation that governs the weight distribution $\pi_t^{(\mu)}$, then deriving the moment equations and equations for mean and covariance, and finally closing the equations via Gaussian filtering.

\subsection{Flow of the Neural Galerkin weight distribution}
Building on a pre-trained CoLoRA network, we now consider random online weight vectors $\Theta(t, \mu)$ that are distributed as  $\Theta(0, \mu) \sim \pi_0^{(\mu)}$ at the initial time $t = 0$. 
Correspondingly, the output $\hat{q}(\Theta(t, \mu), x; \theta_{\text{off}})$ at any $x \in \Omega$ is also a random variable. 
Consider now the Neural Galerkin dynamics given by the velocity field $v$ specified in \eqref{eq:Prelim:NGODE}. 
Recall that the velocity field $v$ induces a flow map $\Phi^{(\mu)}_t: \mathbb{R}^p \to \mathbb{R}^p$. 
The distribution of $\Theta(t, \mu)$ is given by the pushforward of the prior $\pi_0^{(\mu)}$ through the flow $\Phi_t^{(\mu)}$,
\[
\pi_t^{(\mu)} = (\Phi_t^{(\mu)})_{\sharp}\pi_0^{(\mu)}
\]
so that we can write the density $\pi_t^{(\mu)}$ as
\[
\pi_t^{(\mu)}(\theta) = \pi_0^{(\mu)}(\Phi_{-t}^{(\mu)}(\theta)) |\operatorname{det}(\nabla_{\theta} \Phi_{-t}^{(\mu)}(\theta))|\,.
\]
Furthermore, the density $\pi_t^{(\mu)}$ of the weights $\Theta(t, \mu)$ evolves via the continuity equation as
\begin{equation}\label{eq:NGCE}
\partial_t \pi_t^{(\mu)}(\theta) = - \nabla \cdot(v(\theta; \mu) \pi_t^{(\mu)}(\theta))\,,
\end{equation}
with appropriate boundary conditions such as vanishing flux at infinity or compactly supported $\pi_t^{(\mu)}$ so that boundary contributions vanish. We refer to, e.g., \cite{ambrosio_gradient_2005} for details. 
We then have for smooth compactly supported test functions $\varphi: \mathbb{R}^p \to \mathbb{R}$ that
\begin{equation}\label{eq:NG:WeakFormCE}
\frac{\mathrm d}{\mathrm dt} \int \varphi(\theta) \pi_t^{(\mu)}(\theta) \mathrm d\theta = \int \nabla_{\theta} \varphi(\theta) \cdot v(\theta; \mu)\pi_t^{(\mu)}(\theta) \mathrm d\theta
\end{equation}
holds, where we assumed that the function $\pi_t^{(\mu)}$ either has compact support or decays fast enough in the far field so that the boundary terms vanish in \eqref{eq:NG:WeakFormCE}. 
The calculations above that involve the continuity equation \eqref{eq:NGCE} and the weak form \eqref{eq:NG:WeakFormCE} are formal and only meant to motivate the derivation of the moments in the next section. We refer to \cite{DiPerna1989,ambrosio_gradient_2005} for details.

\subsection{Moments of the Neural Galerkin weight distribution}
From the weak form \eqref{eq:NG:WeakFormCE}, we can derive the first and second moment equations by using the standard procedure of using monomials as test functions \cite{StochHandbook,PavliotisBook}.
Let us first consider the first moment 
\begin{equation}\label{eq:FNG:DefFirstM}
m(t, \mu) = \mathbb{E}_{\pi_t^{(\mu)}}[\Theta]\,.
\end{equation}
Testing \eqref{eq:NG:WeakFormCE} with $\varphi_i(\theta) = \theta_i$ leads to 
\begin{equation}\label{eq:FNG:FirstMEq}
\frac{\mathrm d}{\mathrm dt}\mathbb{E}_{\pi_t^{(\mu)}}[\Theta_i] = \mathbb{E}_{\pi_t^{(\mu)}}[e_i \cdot v(\Theta; \mu)]\,,\qquad i = 1, \dots, p\,,
\end{equation}
where $e_i \in \mathbb{R}^p$ denotes the $i$-th unit vector in dimension $p$.
Taking all $p$ equations in \eqref{eq:FNG:FirstMEq} together and using \eqref{eq:FNG:DefFirstM} gives the first moment equation
\begin{equation}\label{eq:NG:MeanEquation}
\frac{\mathrm d}{\mathrm dt} m(t, \mu) = \mathbb{E}_{\pi_t^{(\mu)}}[v(\Theta; \mu)]\,.
\end{equation}
Analogously, we can derive the second moment equations using the test functions 
\[
\varphi_{i,j}(\theta) = \theta_i\theta_j\,,\qquad i,j = 1, \dots, p\,,
\]
to obtain
\[
\frac{\mathrm d}{\mathrm dt}\mathbb{E}_{\pi_t^{(\mu)}}[\Theta_i\Theta_j] = \mathbb{E}_{\pi_t^{(\mu)}}[(\Theta_j e_i + \Theta_ie_j) \cdot v(\Theta; \mu)]\,,\quad i, j = 1, \dots, p\,,
\]
and then the second moment equation in matrix form 
\begin{equation}\label{eq:FNG:2ndMomEq}
\frac{\mathrm d}{\mathrm dt}\mathbb{E}_{\pi_t^{(\mu)}}[\Theta\Theta^{\top}] = \mathbb{E}_{\pi_t^{(\mu)}}[v(\Theta; \mu) \Theta^{\top} + \Theta v(\Theta; \mu)^{\top}]\,.
\end{equation}
The first and the second moment equations can be used to derive the equations for the covariance
\[
\Sigma(t, \mu) = \mathbb{E}_{\pi_t^{(\mu)}}[(\Theta - m(t, \mu))(\Theta - m(t, \mu))^{\top}]
\]
which is
\begin{equation}\label{eq:NG:CovEq}
\frac{\mathrm d}{\mathrm dt}\Sigma(t, \mu) = \mathbb{E}_{\pi_t^{(\mu)}}[v(\Theta; \mu)(\Theta - m(t, \mu))^{\top}] + \mathbb{E}_{\pi_t^{(\mu)}}[(\Theta - m(t, \mu))v(\Theta; \mu)^{\top}]\,.
\end{equation}
The equations for the mean \eqref{eq:NG:MeanEquation} and covariance \eqref{eq:NG:CovEq} are not closed because they depend on the density $\pi_t^{(\mu)}$ which cannot be described by the mean and covariance alone.

\subsection{Filtered Neural Galerkin equations with closed moments}
If the velocity field $v$ is affine in the weight $\theta$, then the moment equations \eqref{eq:NG:MeanEquation}--\eqref{eq:FNG:2ndMomEq} are closed if the initial distribution $\pi_0^{(\mu)}$ is Gaussian; however, the velocity field $v$ is not affine in case of Neural Galerkin schemes because of the nonlinearity of the parametrization $\hat{q}$. 
Instead, we define the filtered Neural Galerkin equations via the first-order approximation of the velocity field, which is a classical approach that can also be used to derive, e.g., the extended Kalman filter \cite{Sarkka_2013,Sanz-Alonso_Stuart_Taeb_2023}. 
We linearize $v$ in $\theta$ at the current mean $m(t, \mu)$ to obtain
\begin{equation}\label{eq:NG:VBar}
\hat{v}(\theta; \mu) = v(m(t, \mu); \mu) + \nabla_{\theta}v(m(t, \mu); \mu)(\theta - m(t, \mu))\,.
\end{equation}
Plugging the first-order approximation $\hat{v}$ into the mean and covariance equations \eqref{eq:NG:MeanEquation} and \eqref{eq:NG:CovEq} gives the filtered Neural Galerkin equations
\begin{equation}\label{eq:FNG:MomentEq}
\begin{aligned}
\frac{\mathrm d}{\mathrm dt}\hat{m}(t, \mu) = & v(\hat{m}(t, \mu); \mu)\,,\\
\frac{\mathrm d}{\mathrm dt}\hat{\Sigma}(t, \mu) = & \nabla_{\theta} v(\hat{m}(t, \mu); \mu)\hat{\Sigma}(t, \mu) + \hat{\Sigma}(t, \mu) \nabla_{\theta} v(\hat{m}(t, \mu); \mu)^{\top}\,,
\end{aligned}
\end{equation}
where $\hat{m}$ and $\hat{\Sigma}$ denote the approximate mean and covariance corresponding to the first-order  dynamics. The initial conditions are $\hat{m}(0, \mu) = \mathbb{E}_{\pi_0^{(\mu)}}[\Theta]$ and $\hat{\Sigma}(0, \mu) = \mathbb{E}_{\pi_0^{(\mu)}}[\Theta\Theta^{\top}] - \mathbb{E}_{\pi_0^{(\mu)}}[\Theta]\mathbb{E}_{\pi_0^{(\mu)}}[\Theta]^{\top}$. System \eqref{eq:FNG:MomentEq} is closed.

We remark that deriving a closed system via first-order dynamics is only one approach for deriving filtered Neural Galerkin dynamics. There are other approaches for closing moment equations which can lead to different variants of filtered Neural Galerkin schemes. We refer to \cite{StochHandbook,PavliotisBook,Sarkka_2013,Sanz-Alonso_Stuart_Taeb_2023} for more extensive discussions of moment equations.

\subsection{Bayesian pre-training of CoLoRA networks}\label{sec:FNG:CoLoRA}
Let us now circle back to pre-training the CoLoRA network. 
The filtered Neural Galerkin dynamics approximate the weight vector as a Gaussian random variable with distribution $\hat{\pi}_t^{(\mu)}$ given by its mean $\hat{m}(t, \mu)$ and covariance matrix $\hat{\Sigma}(t, \mu)$.
Correspondingly, when training a CoLoRA network, the hyper-network $h_{\psi}: (t, \mu) \mapsto [\hat{m}(t, \mu), \hat{\eta}(t,  \mu)]$ maps time $t$ and parameter $\mu$ to the mean and reparametrized components $\hat{\eta}(t, \mu)$ of the covariance matrix. In the following, we focus only on diagonal covariance matrices so that $\hat{\Sigma}_{ii}(t, \mu)= (\log\left(1+\exp\left(\hat{\eta}_i(t,  \mu)\right)\right))^2$ for $i = 1, \dots, p$. 
The offline weights $\theta_{\text{off}}$ remain deterministic. 
To pre-train a CoLoRA network $\hat{q}(\Theta(t, \mu), \cdot; \theta_{\text{off}}): \Omega \to \mathbb{R}$ in preparation for random online weights $\Theta(t, \mu)$, we introduce the likelihood
\begin{equation}\label{eq:FNG:CoLoRALilelihood}
\ell(y\, |\, \theta, x) = \exp\left(- \frac{1}{\sigma^2}\|y - \hat{q}(\theta, x; \theta_{\text{off}})\|_2^2\right)
\end{equation}
with $\sigma > 0$, which imposes a Gaussian noise model on the data. 
The prior of the online weights is $\nu$; it is fixed over all times $t$ and parameters $\mu$. 
Given training data set $\mathcal{D}(t, \mu)$ in \eqref{eq:Dtmu}, the posterior distribution is 
\[
\pi_{\text{post}}(\theta(t, \mu) | \Dcal(t, \mu)) \ \propto\ \nu(\theta(t, \mu))\cdot \prod_{x \in \Omega_{\text{train}}} \ell(q(t, x; \mu) \, | \, \theta(t, \mu), x)\,. 
\]
We then train the CoLoRA network over $\theta_{\text{off}}$ and $\psi$ such that the posterior $\pi_{\text{post}}(\theta(t, \mu) | \Dcal(t, \mu))$ at time $t$ and $\mu$ is close in the Kullback-Leibler divergence to the Gaussian $\pi_{h_{\psi}(t, \mu)}$, where the hyper-network $h_{\psi}$ evaluates via the reparametrization of $(\hat{m}(t, \mu), \hat{\eta}(t, \mu))$ to the mean $\hat{m}(t, \mu)$ and covariance matrix $\hat{\Sigma}(t, \mu)$ of $\pi_{h_{\psi}(t, \mu)}$,
\[
\min_{\theta_{\text{off}} \in \mathbb{R}^n, \psi \in \mathbb{R}^{r'}}\quad \sum_{i = 1}^m\sum_{t \in \{t_0, \dots, t_K\}} \operatorname{KL}(\pi_{h_{\psi}(t, \mu_i)} \,||\,\pi_{\text{post}}(\cdot | \Dcal(t, \mu_i))  )\,,
\]
where the sets $\mathcal{D}(t, \mu_i)$ are defined in \eqref{eq:Dtmu} and $\mu_1, \dots, \mu_m$ are the training parameters. 
It is sufficient to maximize the corresponding evidence lower bound, which is given by
\begin{multline}\label{eq:BColoRAObj}
\mathcal{L}(\theta_{\text{off}}, \psi) = \sum_{i=1}^m\sum_{t \in \{t_0, \dots, t_K\}} \mathbb{E}_{\Theta \sim \pi_{h_{\psi}(t, \mu_i)}} \left[\sum_{x \in \Omega_{\text{train}}} \log \ell(q(t, x; \mu_i) \,|\, \Theta, x)\right]\\ - \mathbb{E}_{\Theta \sim \pi_{h_{\psi}(t, \mu_i)}}[\log \pi_{h_{\psi}(t, \mu_i)}(\Theta) - \log \nu(\Theta)]\,.
\end{multline}
Once the CoLoRA network is pre-trained with the objective \eqref{eq:BColoRAObj}, the hyper-network $h_{\psi}$ can be used to generate the mean $\hat{m}(0, \mu)$ and covariance $\hat{\Sigma}(0, \mu)$ at time $t = 0$ to define the distribution $\hat{\pi}_0^{(\mu)}$ of the weights $\Theta(0, \mu)$.  
We remark that the hyper-network $h_{\psi}$ can also be used to generate mean and covariance at later times $t > 0$, which leads to a purely data-driven approach to push forward the distribution of the weights. We do not  pursue this direction further in this work but refer to \cite{pmlr-v235-berman24b} that considers this data-driven perspective in the setting where the initial condition is deterministic.

\subsection{Computational procedure of filtered Neural Galerkin schemes}
We now  derive a computational procedure for filtered Neural Galerkin schemes with pre-trained CoLoRA networks. 

\subsubsection{Computational procedure}
In the offline phase, the CoLoRA network is trained with the procedure described in \Cref{sec:FNG:CoLoRA} to obtain the offline weights $\theta_{\text{off}}$ and the hyper-network $h_{\psi}$. In the online phase, for a new parameter $\mu \in \mathcal{Q}$,  
the filtered Neural Galerkin equations \eqref{eq:FNG:MomentEq} are solved to approximate $\hat{\pi}_t^{(\mu)}$ at discrete time points in the time interval $[0, T]$. For demonstration purposes, we use here explicit Euler time integration to discretize the filtered Neural Galerkin equations, though we use Runge-Kutta 4 as time integration scheme in our numerical experiments. Let $\delta t > 0$ be the time-step size, then we obtain
\begin{equation}\label{eq:FNG:DiscFNG}
\begin{aligned}
\hat{m}_{k + 1}(\mu) = & \hat{m}_k(\mu) + \delta t v(\hat{m}_k(\mu); \mu)\,,\\
\hat{\Sigma}_{k + 1}(\mu) = &\hat{\Sigma}_k(\mu) + \delta t\left(\nabla_{\theta}v(\hat{m}_k(\mu); \mu)\hat{\Sigma}_k(\mu) + \hat{\Sigma}_k(\mu)\nabla_{\theta}v(\hat{m}_k(\mu); \mu)^{\top}\right)\,,
\end{aligned}
\end{equation}
for time steps $k = 0, \dots, K - 1$ and $K = \lceil T/\delta t\rceil$. The initial mean $\hat{m}_0(\mu)$ and covariance $\hat{\Sigma}_0(\mu)$ are computed with the hyper-network $h_{\psi}$, where the hyper-network outputs the reparametrized covariance; see \Cref{sec:FNG:CoLoRA}. Once the mean $\hat{m}_k(\mu)$ and covariances $\hat{\Sigma}_k(\mu)$ have been computed over all time steps $k = 0, \dots, K$, they define the Gaussian distribution $\hat{\pi}_k^{(\mu)}$ at the steps $k = 0, \dots, K$. 

Evaluating the right-hand side of \eqref{eq:FNG:DiscFNG} requires evaluating the velocity field $v$. To evaluate the velocity field $v$, we select a set $\{x_1, \dots, x_N\} \subset \Omega$ of $N \in \mathbb{N}$ collocation points to form the batch gradient
\begin{equation}\label{eq:FNG:BatchGradient}
J(\theta) = \begin{bmatrix}
- & \nabla_{\theta}\hat{q}(\theta, x_1; \theta_{\text{off}})^{\top} & -\\
& \vdots &\\
- & \nabla_{\theta}\hat{q}(\theta, x_N; \theta_{\text{off}})^{\top} & -
\end{bmatrix} \in \mathbb{R}^{N \times p}
\end{equation}
and the batch right-hand side 
\begin{equation}\label{eq:FNG:BatchRHS}
\bar{f}(\theta; \mu) = \begin{bmatrix}
f(x_1, \hat{q}(\theta, \cdot; \theta_{\text{off}}); \mu) & \dots & f(x_N, \hat{q}(\theta, \cdot; \theta_{\text{off}}); \mu) \end{bmatrix}^{\top} \in \mathbb{R}^N
\end{equation}
and compute the value of $v(\theta; \mu)$ via the least-squares problem $\min_v \|J(\theta) v - \bar{f}(\theta; \mu)\|_2$. Thus, the velocity field is implicitly defined via the solution of the least-squares problem. Differentiating the velocity field with respect to $\theta$ to obtain $\nabla_{\theta} v$ can be achieved by differentiating through the least-squares solution \cite{NEURIPS2022_228b9279}, which is implemented by common automatic differentiation tools such as JAX \cite{jax2018github}.  

\subsubsection{Cost complexity}
Let us first consider the cost complexity of a time step with the Neural Galerkin scheme following the dynamics \eqref{eq:Prelim:VField}. The batch gradient \eqref{eq:FNG:BatchGradient} and the batch right-hand side \eqref{eq:FNG:BatchRHS} are formed and the corresponding least-squares problem to evaluate $v$ is computed, which incurs costs that scale as $\mathcal{O}(Np^2)$, where $p$ is the dimension of the online weight vector. So if an ensemble of $M$ members is computed, this leads to a cost scaling of $\mathcal{O}(MNp^2)$.

When solving the filtered equations \eqref{eq:FNG:DiscFNG} at time step $k$, we have to compute $v$ at $\hat{m}_k(\mu)$ as well as apply $\nabla_{\theta}v$ to $\hat{\Sigma}_k(\mu)$, which incurs costs $\mathcal{O}(Np^2)$ for the least-squares solve to obtain the value of $v$ and $\mathcal{O}(Np^3)$ for applying $\nabla_{\theta}v$: For each column of $\hat{\Sigma}_k(\mu)$, we have to multiply with $\nabla_{\theta}v$. For each multiplication, we have to differentiate in the corresponding direction through the least-squares problem underlying the evaluation of $v$, which scales as $\mathcal{O}(Np^2)$. Thus, because we repeat this for all of the $p$ columns of $\hat{\Sigma}_k(\mu)$, we obtain costs that scale as $\mathcal{O}(Np^3)$; see, e.g., \cite{NEURIPS2022_228b9279} for details on implicit differentiation of least-squares problems. Computing $\hat{m}_{k + 1}(\mu)$ at the next time step $k + 1$ incurs costs that scale as $\mathcal{O}(p)$ but the costs of computing $\hat{\Sigma}_{k + 1}(\mu)$ scale as $\mathcal{O}(p^3)$, because matrix-matrix multiplications with matrices of size $p \times p$ are required. Thus, the costs of one time step with filtered Neural Galerkin scale as $\mathcal{O}(Np^2 + Np^3 + p^3)$. Because typically $p \ll N$, recall that $p$ is the dimension of the online weights, we obtain that costs scale as $\mathcal{O}(Np^3)$. Thus, with filtered Neural Galerkin schemes we incur an extra factor $p$ in the cost complexity compared to one Neural Galerkin time step. However, recall that $p$ is the dimension of the online weights, which typically can be chosen lower than the number of ensemble members.

\section{Numerical experiments}\label{sec:NumExp}

We demonstrate filtered Neural Galerkin schemes on two problems:
the one-dimensional viscous Burgers equation and a particle problem governed by the two-dimensional Vlasov equation.

\subsection{Viscous Burgers equation}
We consider the viscous Burgers equation. 
\subsubsection{Burgers: Setup}
\label{sec:burgers}
Let us consider the one-dimensional viscous Burgers equation on a periodic spatial domain,
\begin{equation}
\partial_t q(t,x;\mu) 
+ \tfrac{1}{2}\partial_x \big(q(t,x;\mu)^2\big) 
= \mu\,\partial_{xx} q(t,x;\mu),
\qquad x \in [-1,1), \quad t \in [0,T],
\label{eq:burgers-pde}
\end{equation}
with viscosity $\mu > 0$ and end time $T>0$. The initial condition is 
\begin{equation}
q(0,x;\mu)
= 0.8 \exp\left(-\frac{(x-x_1)^2}{\sigma^2}\right),
\qquad x_1 = -0.2,\quad \sigma = 0.2,
\label{eq:burgers-ic}
\end{equation}
which is numerically zero at the boundaries. We discretize the PDE with a second-order finite-volume scheme in space on a grid with $256$ grid points. In time, we use a Runge-Kutta 4 with time-step size in $\delta t \in [10^{-4}, 3 \times 10^{-3}]$ depending on the viscosity parameter $\mu$. 
For each viscosity $\mu$, we solve until end time $T=2.0$. We generate training data \eqref{eq:Prelim:TrainingTraj} for $m = 34$ equidistant $\mu_1, \dots, \mu_m$ in the parameter domain $[10^{-3},\,10^{-1}]$. The set $\Omega_{\text{train}}$ contains the 256 equidistant grid points in the spatial domain $\Omega$ and in time we consider the equidistant time points $0 = t_0 < t_1 < \dots < t_K = 2.0$ with $K = 100$. Altogether, we have sets \eqref{eq:Dtmu} for $(t, \mu)$ pairs in $\{t_0, t_1, \dots, t_k\} \times \{\mu_1, \dots, \mu_m\}$. 

\begin{figure}[t]
\centering\begin{tabular}{cc}
\includegraphics[width=0.45\columnwidth]{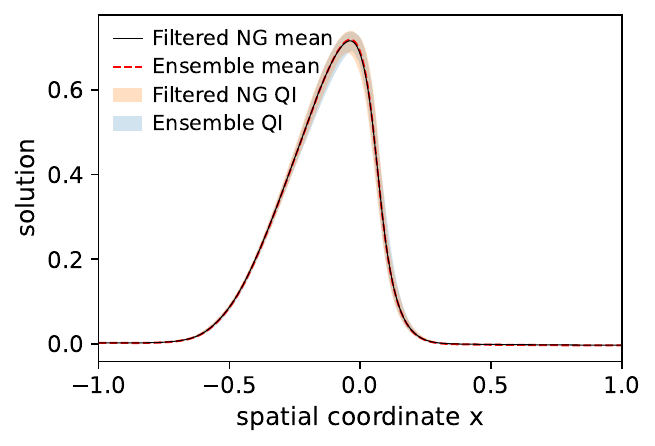} & \includegraphics[width=0.45\columnwidth]{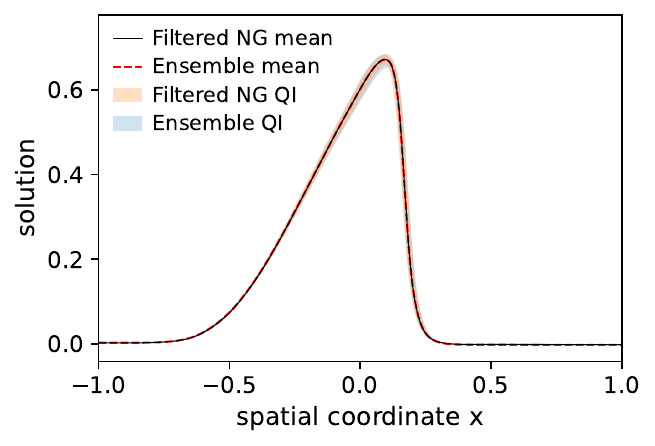}\\
\scriptsize (a) $t = 0.25$ & \scriptsize (b) $t = 0.5$\\
\includegraphics[width=0.45\columnwidth]{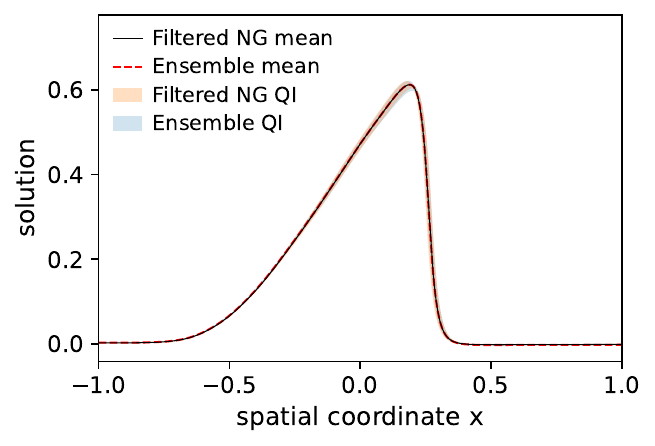} & \includegraphics[width=0.45\columnwidth]{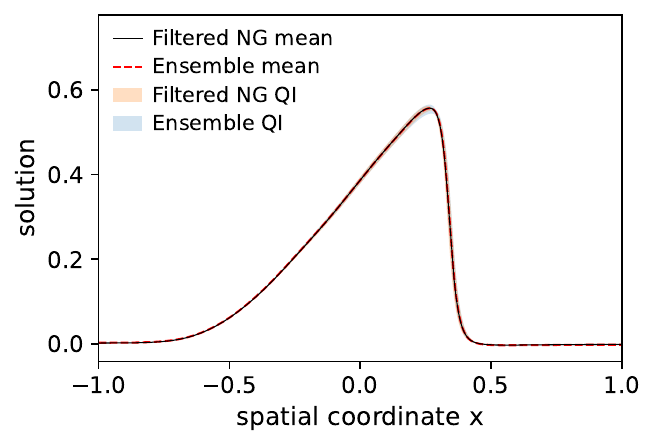}\\
\scriptsize (c) $t = 0.75$ & \scriptsize (d) $t = 1.0$
\end{tabular}
\caption{Burgers: A solution of the filtered Neural Galerkin (NG) model provides comparable quantile intervals (QI) as an ensemble of 100 Neural Galerkin solutions for different realizations of the initial condition in this example. }
\label{fig:Burgers:ResultSmallBand}
\end{figure}

\begin{figure}[t]
\centering\begin{tabular}{cc}
\includegraphics[width=0.45\columnwidth]{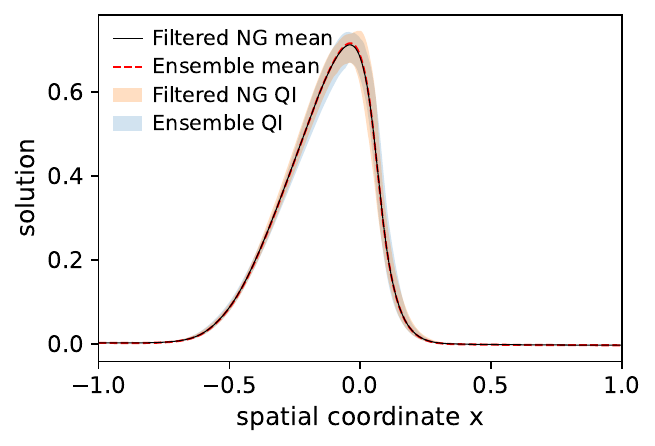} & \includegraphics[width=0.45\columnwidth]{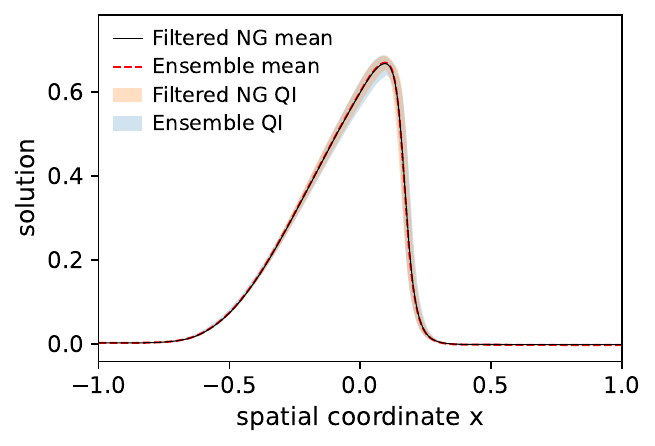}\\
\scriptsize (a) $t = 0.25$ & \scriptsize (b) $t = 0.5$\\
\includegraphics[width=0.45\columnwidth]{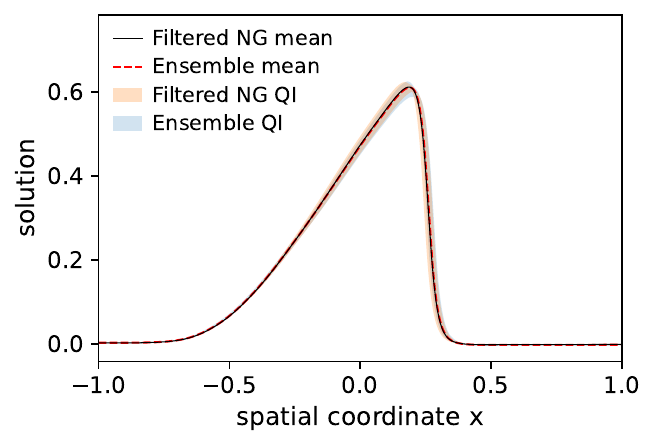} & \includegraphics[width=0.45\columnwidth]{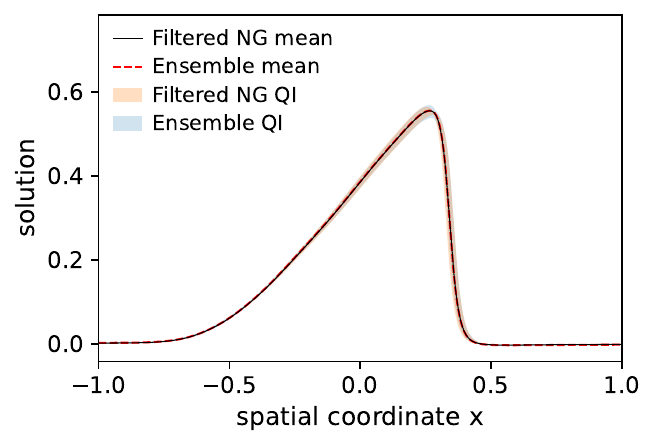}\\
\scriptsize (c) $t = 0.75$ & \scriptsize (d) $t = 1.0$
\end{tabular}
\caption{Burgers: The filtered Neural Galerkin approach approximates well the mean and quantile interval (QI) that is generated by an ensemble of 100 Neural Galerkin solutions.}
\label{fig:Burgers:ResultWideBand}
\end{figure}

\subsubsection{Burgers: Bayesian pre-training of CoLoRA networks}\label{sec:NumExp:Burgers:Pretrain}
We train a CoLoRA network as described in \Cref{sec:FNG:CoLoRA}. The number of hidden layers is $L = 5$, the width is $64$ and the activation function is $\tanh$. The output layer is given in \eqref{eq:Prelim:CoLoRAOutputLayer}. Overall, the online weight vector has dimension $p = L+1= 6$. 
For the hidden layers, we set the CoLoRA rank to $r = 16$. 
The hyper-network is a fully connected feedforward network with four layers, each with width $64$. Hidden layer number one and three have ReLU activation functions and hidden layer number two has the sigmoid function as activation function. The output layer is linear. As collocation points we use the $N = 1000$ equidistant points in $\Omega$. 
For training with the loss \eqref{eq:BColoRAObj}, we replace the expectations with Monte Carlo estimators with 
ten samples. The prior $\nu$ is a standard normal and the noise standard deviation is $\sigma_{\text{noise}}=0.01$ in the likelihood \eqref{eq:FNG:CoLoRALilelihood}.
We train for $300$ epochs using the Adam optimizer with learning rate $10^{-3}$ and batch size $4096$.

\subsubsection{Burgers: Results}\label{sec:NumExp:Burgers:Results}
We now compare the 95\% quantile interval computed with filtered Neural Galerkin to the interval obtained from ensembles of Neural Galerkin solutions. In filtered Neural Galerkin, we approximate $\pi_t^{(\mu)}$ with a Gaussian $\hat{\pi}_t^{(\mu)}$ with mean $\hat{m}(t, \mu)$ and covariance $\hat{\Sigma}(t, \mu)$ over time $t$. To empirically estimate the 95\% quantile interval at time $t$, we draw 100 samples from $\hat{\pi}^{(\mu)}_t$ and evaluate the corresponding pre-trained CoLoRA network. We keep 95\% of the 100 samples and show in the following the corresponding interval that these samples span. If the matrix $\hat{\Sigma}(t, \mu)$ is not symmetric positive definite, we symmetrize it and truncate all non-positive eigenvalues before drawing samples. 
We compare to ensembles of Neural Galerkin solutions. To generate an ensemble member, we draw a realization of the initial condition $\pi_0^{(\mu)}$ and integrate the Neural Galerkin equations \eqref{eq:Prelim:NGODE} with the corresponding pre-trained CoLoRA network. We generate 100 ensemble members and compute the empirical 95\% quantile interval. 

Let us consider the test parameter $\mu^{\text{test}} = 0.005$. \Cref{fig:Burgers:ResultSmallBand} shows the quantile intervals that are obtained at times $t \in \{0.25, 0.5, 0.75, 1.0\}$. The quantile intervals obtained with filtered Neural Galerkin are comparable to the ones obtained with ensembles of Neural Galerkin solutions. For obtaining the result shown in \Cref{fig:Burgers:ResultWideBand}, we increase the noise in the likelihood \eqref{eq:FNG:CoLoRALilelihood} to $\sigma_{\text{noise}} = 0.015$ to have an initial distribution $\hat{\pi}_0^{(\mu)}$ with a larger variance. Correspondingly, the larger variance is propagated forward. Our filtered Neural Galerkin schemes compute comparable quantile intervals as obtained with ensembles of solutions. 

\begin{figure}[t]
\centering\begin{tabular}{cc}
\includegraphics[width=0.45\columnwidth]{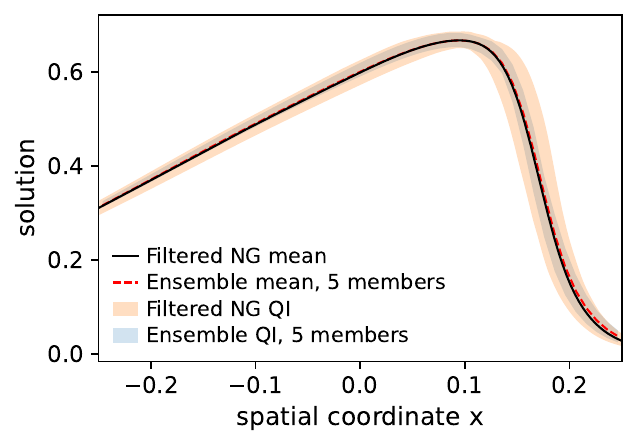} & \includegraphics[width=0.45\columnwidth]{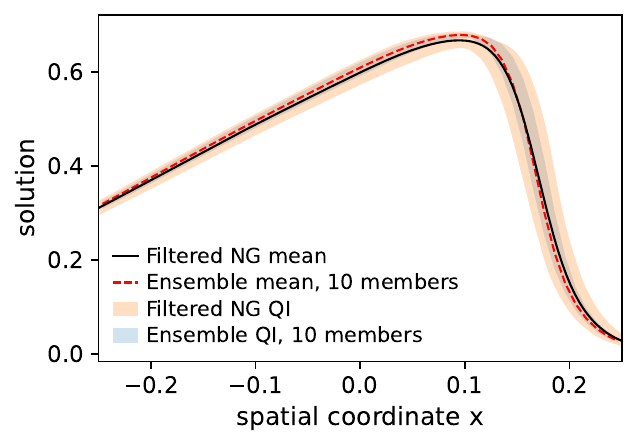}\\
\scriptsize (a) 5 ensemble members & \scriptsize (b) 10 ensemble members\\
\includegraphics[width=0.45\columnwidth]{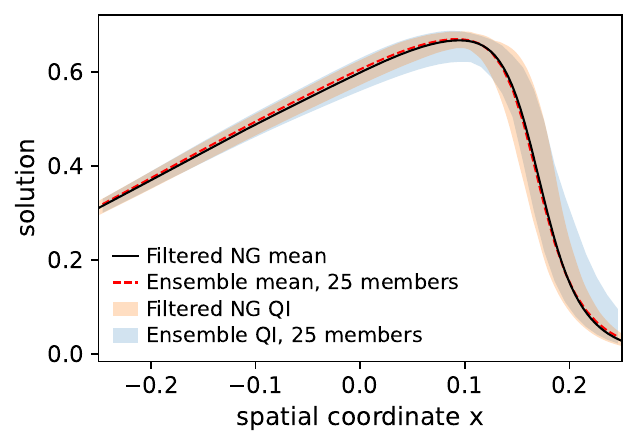} & \includegraphics[width=0.45\columnwidth]{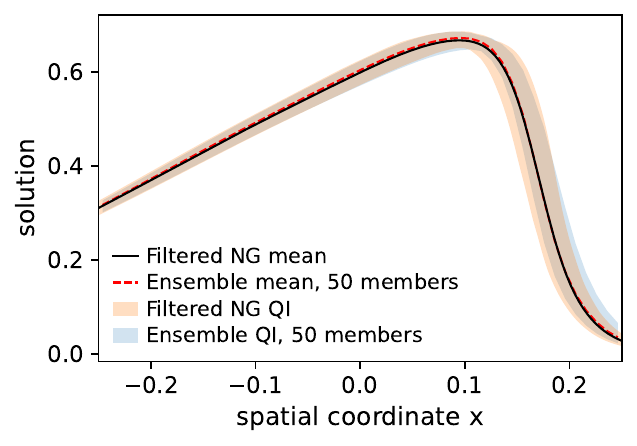}\\
\scriptsize (c) 25 ensemble members & \scriptsize (d) 50 ensemble members
\end{tabular}
\caption{Burgers: At least 50 ensemble members are needed to achieve a comparable accuracy as the quantile interval (QI) predicted by filtered Neural Galerkin schemes. Time is $t =0.5$.}
\label{fig:Burgers:EnsembleMembers}
\end{figure}

\begin{figure}[t]
\centering\begin{tabular}{cc}
\includegraphics[width=0.45\columnwidth]{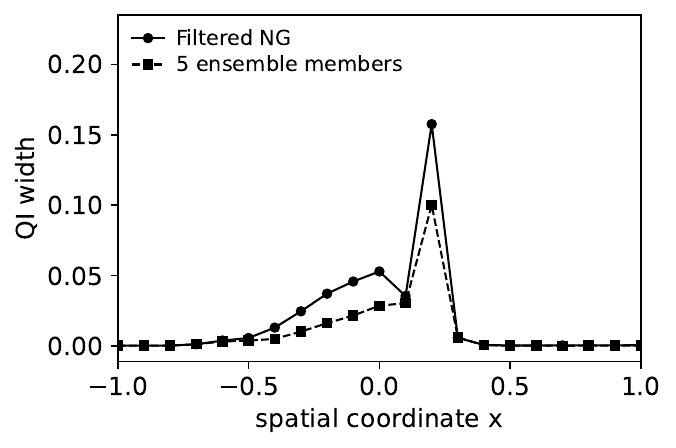} & 
\includegraphics[width=0.45\columnwidth]{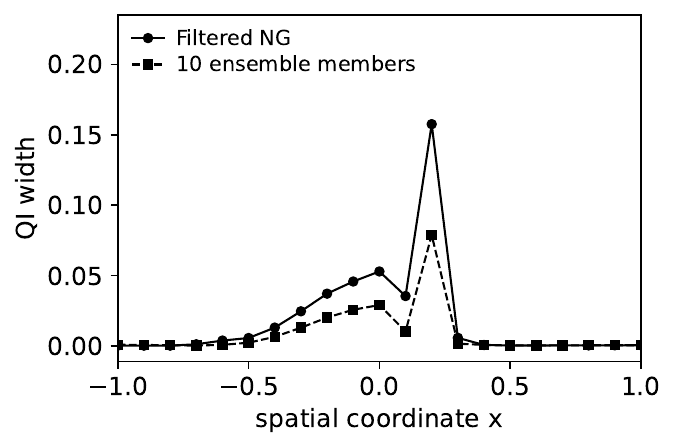}\\
\scriptsize (a) 5 ensemble members & \scriptsize (b) 10 ensemble members\\
\includegraphics[width=0.45\columnwidth]{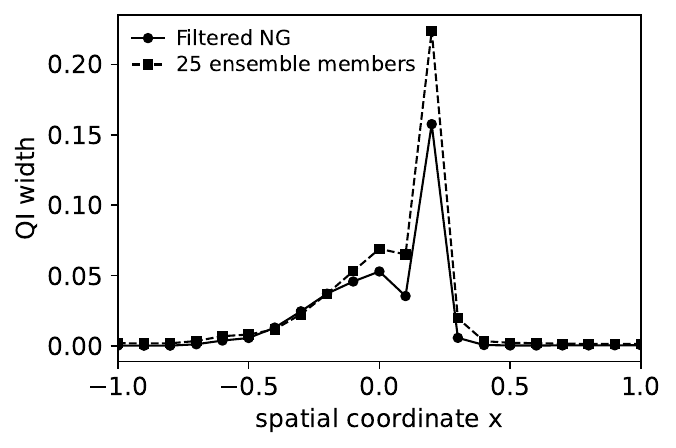} &
\includegraphics[width=0.45\columnwidth]{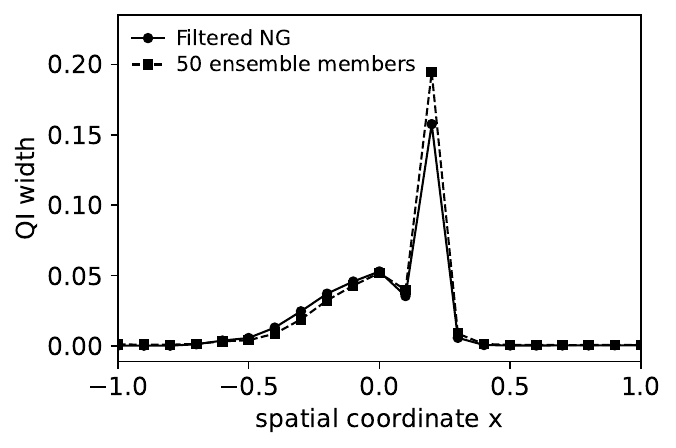}\\
\scriptsize (c) 25 ensemble members &
\scriptsize (d) 50 ensemble members
\end{tabular}
\caption{Burgers: Shows that the width of the quantile intervals computed with filtered Neural Galerkin and an ensemble of Neural Galerkin solutions starts to agree when at least 50 ensemble members are used. }
\label{fig:Burgers:EnsembleMembersWidth}
\end{figure}

\Cref{fig:Burgers:EnsembleMembers} shows that the quantile interval obtained from 25 ensemble members is still inaccurate compared to the one obtained from 50 members, which indicates that about 50 ensemble members are needed to obtain an accurate interval in this example. In contrast, using only five or ten ensemble members leads to poor predictions of the quantile interval. Additional evidence is provided in \Cref{fig:Burgers:EnsembleMembersWidth} that shows that the width of the quantile intervals computed with filtered Neural Galerkin and an ensemble starts to agree when at least 50 ensemble members are used.  

\begin{figure}
\centering
\begin{tabular}{cc}
\includegraphics[width=0.45\columnwidth]{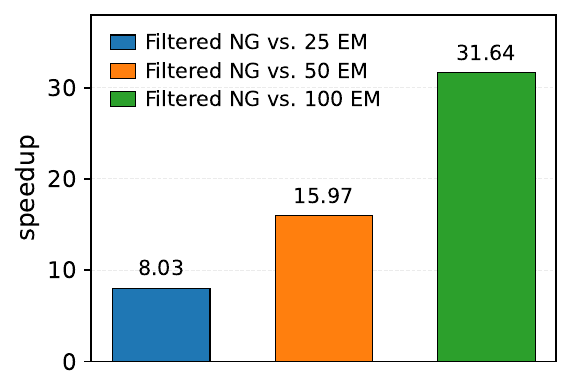} &
\includegraphics[width=0.45\columnwidth]{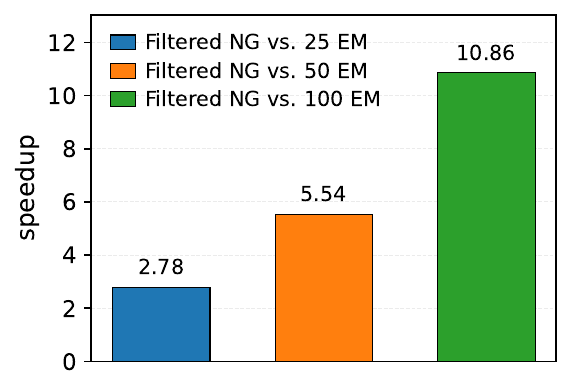}\\
\scriptsize (a) Burgers equation & \scriptsize (b) Vlasov equation
\end{tabular}

\caption{Burgers: Filtered Neural Galerkin generates quantile intervals with more than one order of magnitude lower runtime compared to generating quantile intervals with 100 ensembles members (EM) of Neural Galerkin solutions.}
\label{fig:Burgers:ResultRuntime}
\end{figure}

Taking into account that about 50 ensemble members are needed in this example to obtain an accurate quantile interval, we show in \Cref{fig:Burgers:ResultRuntime} the runtime speedup obtained with filtered Neural Galerkin over ensemble runs with 25, 50, and 100 ensemble members. The filtered Neural Galerkin achieves a speedup of more than one order of magnitude compared to 50 ensemble members in this example.

    \begin{figure}[t]
\begin{tabular}{c}
\includegraphics[width=1.0\columnwidth]{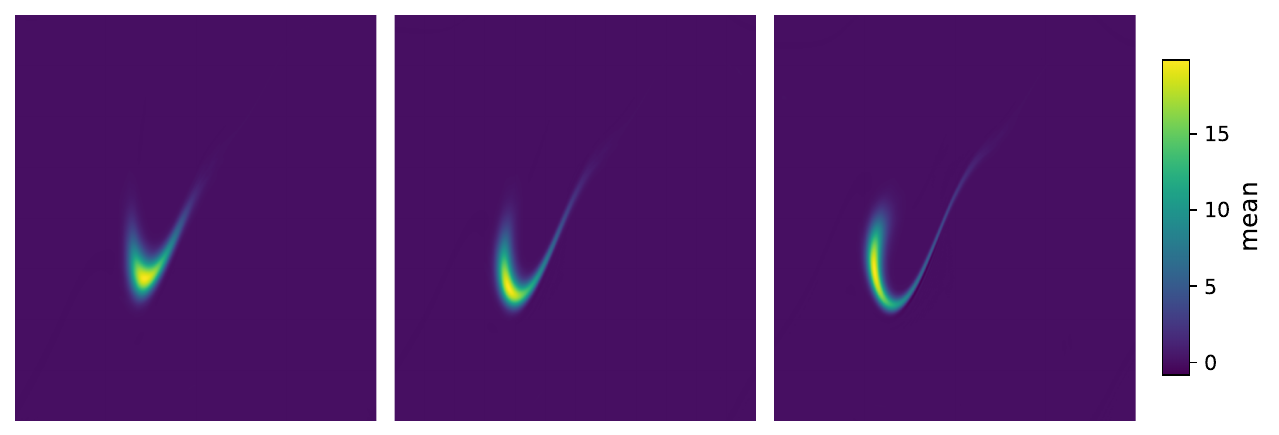}\\
\scriptsize (a) ensemble mean of Neural Galerkin solutions at times $0.50$ (left), $0.75$ (middle), $1.00$ (right)\\
\includegraphics[width=1.0\columnwidth]{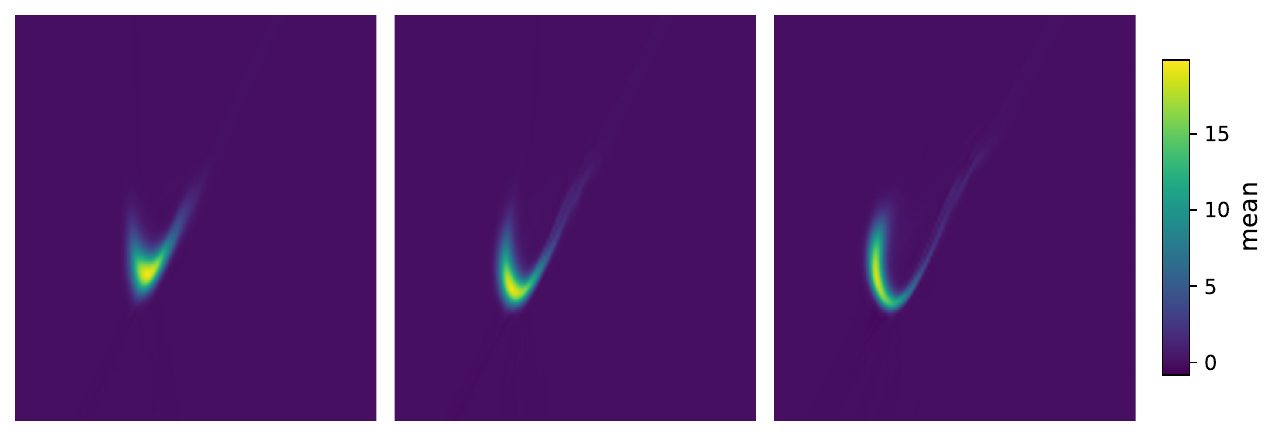}\\
\scriptsize (b) mean solution obtained with filtered Neural Galerkin at times $0.50$ (left), $0.75$ (middle), $1.00$ (right)
\end{tabular}
\caption{Vlasov: The mean solution obtained with filtered Neural Galerkin schemes is in close agreement with the mean of an ensemble of solutions.}
\label{fig:NumExp:Vlasov:ComparisonMean}
\end{figure}

\subsection{Two-dimensional Vlasov equation}
\label{sec:vlasov}
We now consider particle systems that are governed by the two-dimensional Vlasov equation. 
\subsubsection{Vlasov: Setup and Bayesian pre-training of CoLoRA networks}
We follow the setup described in \cite{WEDER2025114249}. The Vlasov equation is
\[
\partial_t q(t,x,v;\mu)
+ v\,\partial_x q(t,x,v;\mu)
- \phi_x(x;\mu)\,\partial_v q(t,x,v;\mu)
= 0,
\qquad (x,v)\in[-1,1)\times[-1,1),
\]
with periodic boundary conditions in both $x$ and $v$. The external field is 
\[
\phi_x(x;\mu)
=
4\alpha \pi \sin\left(\pi(x+\mu)\right)\,\cos^3\left(\pi(x+\mu)\right)
\;-\;
\beta \pi \cos(\pi x),
\]
with fixed $\alpha=0.2$ and $\beta=0.1$. The initial condition is,
\[
q(0,x,v;\mu)
=
\frac{1}{2\pi\gamma}
\exp\left(
-\frac{1}{\pi\gamma}
\left[
\sin^2\left(\frac{\pi}{2}(x-x_0)\right)
+
\sin^2\left(\frac{\pi}{2}(v-v_0)\right)
\right]
\right),
\]
with $x_0=-0.2$, $v_0=0$, and $\gamma=8\times 10^{-3}$. The spatial domain is discretized on an equidistant grid of size $512 \times 512$. Derivatives are discretized with fourth-order finite-difference stencils. The time-step size is $\delta t = 10^{-4}$ and the time integration scheme is Runge-Kutta 4. For generating the training data, we consider the $m = 5$ equidistant parameters $\mu_1, \dots, \mu_m$ in $[0.25, 0.45]$. We set end time to $T = 1$ and the discrete time points at which we store training snapshots to the $K = 50$ equidistant points $0 = t_0 < t_1 < \dots < t_K = T$. 
For training the CoLoRA network, we consider $\Omega_{\text{train}}$ with points corresponding to the $128 \times 128$ equidistant grid in $\Omega$, which are also the collocation points used in Neural Galerkin in this example. The CoLoRA network architecture and pre-training setup is the same as the one described in \Cref{sec:NumExp:Burgers:Pretrain}, except that the first layer now maps from $\mathbb{R}^2$ instead of $\mathbb{R}$ to $\mathbb{R}^{64}$ to account for $x$ and $v$ and the batch size is increased to $8192$.

\subsubsection{Vlasov: Results}
\begin{figure}[t]
\begin{tabular}{c}
\includegraphics[width=1.0\columnwidth]{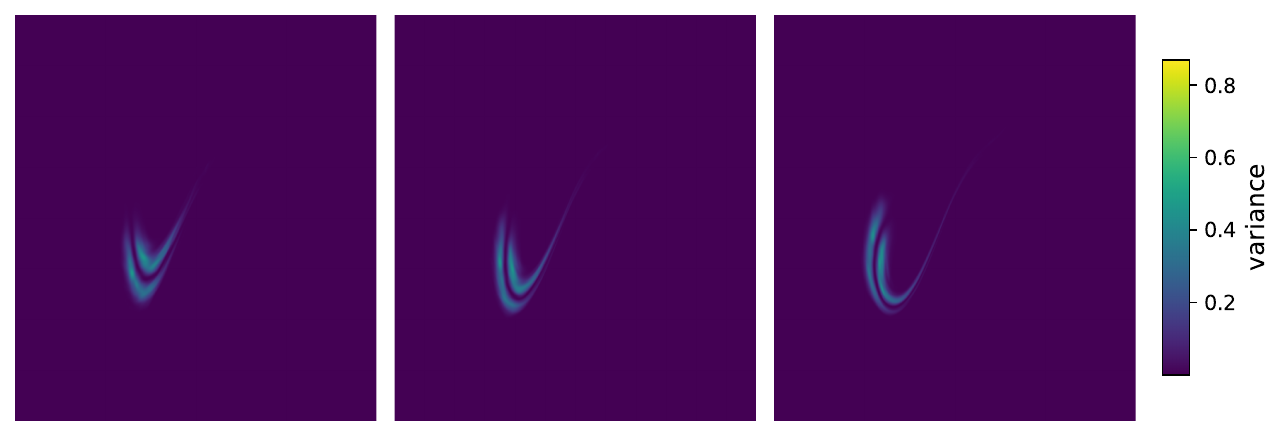}\\
\scriptsize (a) variance of ensemble of Neural Galerkin solutions at times $0.50$ (left), $0.75$ (middle), $1.00$ (right)\\
\includegraphics[width=1.0\columnwidth]{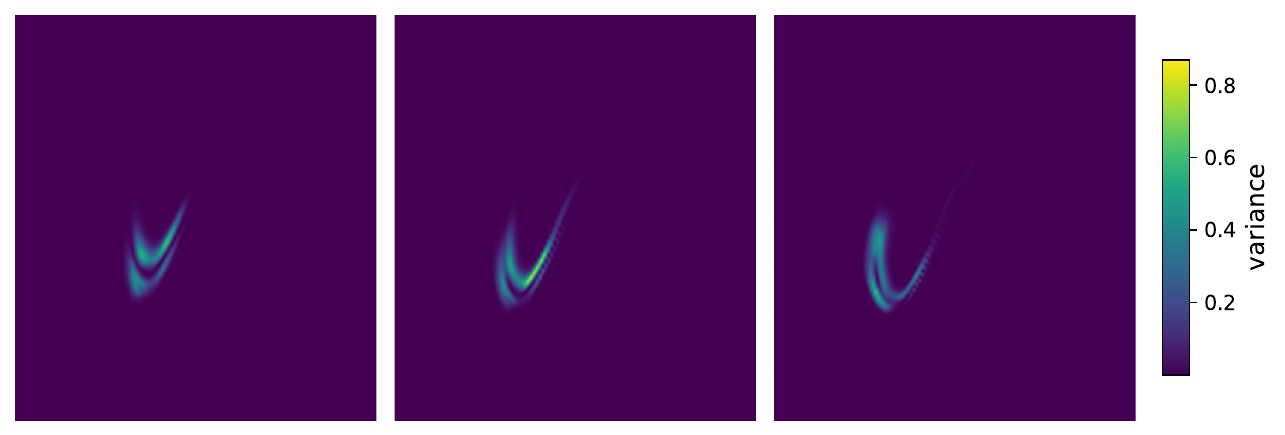}\\
\scriptsize (b) variance obtained with filtered Neural Galerkin at times $0.50$ (left), $0.75$ (middle), $1.00$ (right)
\end{tabular}
\caption{Vlasov: Filtered Neural Galerkin schemes provide variance estimates that are in close agreement to estimates from ensembles of solutions.}
\label{fig:NumExp:Vlasov:ComparisonStd}
\end{figure}

We use filtered Neural Galerkin and ensemble runs of Neural Galerkin solutions to compute a mean solution and its variance. For filtered Neural Galerkin, analogously to \Cref{sec:NumExp:Burgers:Results}, we draw 100 samples from $\hat{\pi}_t^{(\mu^{\text{test}})}$ for the test parameter $\mu^{\text{test}} = 0.375$ and then evaluate the CoLoRA network at the 100 samples to compute mean and variance. Similarly, for an ensemble of 100 Neural Galerkin solutions with initial condition weight vectors sampled from $\pi_0^{(\mu^{\text{test}})}$, we compute mean and variance. The mean solutions are plotted in \Cref{fig:NumExp:Vlasov:ComparisonMean} and the variances are plotted in \Cref{fig:NumExp:Vlasov:ComparisonStd}. The mean and variance obtained with filtered Neural Galerkin is in agreement with the mean and variance obtained from the ensemble. At the same time, the filtered Neural Galerkin achieves a speedup of up to one order of magnitude compared to the ensemble run with 100 samples; see \Cref{fig:Burgers:ResultRuntime}.

\section{Conclusions}\label{sec:Conc}
Filtered Neural Galerkin schemes offer an alternative to ensemble methods for propagating initial condition uncertainties.
We  showed that propagating uncertainties can be achieved with lower costs than typically incurred by ensemble methods, which can be prohibitively expensive even when using reduced models in the context of digital twins. 
Filtered Neural Galerkin schemes avoid ensembles by propagating closed moment equations to obtain a Gaussian approximation of the uncertainties. In numerical experiments, we obtained speedups of at least one order of magnitude compared to ensemble methods.

\section*{Acknowledgments} The authors have been partially funded by the Air Force Office of Scientific Research (AFOSR), USA, award FA9550-24-1-0327.

\bibliographystyle{plain}
\bibliography{main}

\end{document}